\documentclass[11pt,a4paper]{article}
\usepackage{amssymb,amsmath,amsfonts}
\usepackage{url}
\usepackage{epsfig}
\usepackage{hyperref}
\usepackage{srcltx}
\usepackage{xcolor}

\def\poly{{\mathrm{poly}}}

\def\tLin{\mathrm{lin}}
\newtheorem{lemma}{Lemma}[section]
\newtheorem{proposition}{Proposition}[section]
\newtheorem{corollary}{Corollary}[section]

\newtheorem{example}{Example}[section]

\def\RiskOpt{{\hbox{\rm RiskOpt}}}

\def\lin{{\hbox{\rm\tiny lin}}}

\def\three?{3}
\def\four?{4}
\def\ten?{10}

\def\Argmin{\mathop{\hbox{\rm Argmin}}}
\def\beq{\begin{equation}}
\def\eeq{\end{equation}}

\def\norm2to2{{\|\cdot\|_{2,2}}}
\def\Prob{\hbox{\rm Prob}}

\def\bE{{\mathbf{E}}}

\def\Diag{\hbox{\rm  Diag}}
\def\Prob{\hbox{\rm  Prob}}

\def\Opt{\hbox{\rm Opt}}

\def\Tr{{\mathop{\hbox{\rm  Tr}}}}

\def\cB{{\cal B}}

\def\cF{{\cal F}}
\def\cG{{\cal G}}

\def\cN{{\cal N}}
\def\cO{{\cal O}}

\def\cS{{\cal S}}
\def\cT{{\cal T}}

\def\cX{{\cal X}}

\def\Argmin{\mathop{\hbox{\rm  Argmin}}}

\def\Det{{\mathop{\hbox{\rm  Det}}}}

\def\bS{{\mathbf{S}}}

\def\qed{\ \hfill$\square$\par\smallskip}

\def\mypict3{\epsfxsize=220pt\epsfysize=80pt\epsfbox}

\def\bR{{\mathbf{R}}}

\def\argmin{\mathop{\hbox{\rm argmin}}}

\def\Col{{\hbox{\rm Col}}}

\def\Risk{{\hbox{\rm Risk}}}
\newcommand{\hide}[1]{{}}

\newcommand{\be}{\begin{eqnarray}}
\newcommand{\ee}[1]{\label{#1}\end{eqnarray}}
\newcommand{\nn}{\nonumber \\}
\newcommand{\ese}{\end{eqnarray*}}
\newcommand{\bse}{\begin{eqnarray*}}
\newcommand{\rf}[1]{~(\ref{#1})}
\def\SG{{\cal SG}}

\def\half{\tfrac{1}{2}}
\def\four{\tfrac{1}{4}}
\newcommand{\wh}[1]{{\widehat{#1}}}

\def\mr{{{\mathfrak{r}}}}
\def\mp{{{\mathfrak{p}}}}

\def\wt#1{\widetilde{#1}}
\def\ov#1{\overline{#1}}
\def\ul#1{\underline{#1}}
\begin{document}
\title{First order algorithms for computing  linear and polyhedral estimates}
\author{Yannis Bekri,
 LJK, Universit\'e Grenoble Alpes\\ Campus de Saint-Martin-d'H\`{e}res, 38401 France\\
              yannis.bekri@univ-grenoble-alpes.fr
\and Anatoli Juditsky, LJK, Universit\'e Grenoble Alpes\\ Campus de Saint-Martin-d'H\`{e}res, 38401 France\\
              anatoli.juditsky@univ-grenoble-alpes.fr
\and Arkadi Nemirovski, Georgia Institute of Technology\\ Atlanta, Georgia
30332, USA\\
nemirovs@isye.gatech.edu}
\maketitle
\begin{abstract}
It was recently shown \cite{JudNem2018,juditsky2020polyhedral} that ``properly built'' linear and polyhedral estimates  nearly attain minimax accuracy bounds in the problem of recovery of unknown signal  from noisy observations of linear images of the signal when the signal set is an {\em ellitope}. However, design of nearly optimal estimates relies upon solving semidefinite optimization problems with matrix variables, what puts the synthesis of such estimates beyond the rich of the standard Interior Point algorithms of semidefinite  optimization even for moderate size recovery problems. Our goal is to develop First Order Optimization algorithms for the computationally efficient design of linear and polyhedral estimates. In this paper we (a) explain how to eliminate matrix variables, thus reducing dramatically the design dimension when passing from Interior Point to First Order optimization algorithms and  (2) develop and analyse a dedicated algorithm of the latter type---Composite
Truncated Level method.
\end{abstract}


\section{Introduction}
In this paper we discuss numerical algorithms for construction of ``presumably good'' estimates in linear inverse problems.
Specifically, consider the estimation problem as follows. We are given an observation $\omega\in \bR^m$,
\beq\label{eq1}
\omega=Ax+\xi\eeq
where
$\xi\in\bR^m$ is zero mean {random noise},
$A\in \bR^{m\times n}$ is observation matrix, and $x$ is unknown signal known to belong to a  given convex set $\cX\subset \bR^n$. Our objective is to recover the linear image $w=Bx$, $B\in \bR^{\nu\times n}$, of $x$.
\par
Our focus is on {\em linear} and {\em polyhedral estimates} for solving the problem in question.

When applied to the estimation problem above, {\em linear estimate} $\wh w^H_\lin(\omega)$ of $w$ is of the form $ \wh w^H_\lin(\omega)=H^T\omega$ where {\em contrast matrix} $H\in\bR^{m\times \nu}$ is the estimate's parameter. A polyhedral estimate $\widehat{w}^H_\poly(\omega)$ is specified by a {\sl contrast matrix} $H\in\bR^{m\times M}$ according to
\[
\omega\mapsto\wh x^H(\omega)\in\Argmin_{x\in\cX} \left\{\|H^T(\omega-Ax)\|_\infty\right\},\;\; \widehat{w}^H_\poly(\omega) :=B\wh{x}(\omega).
\]
Our interest in these two types of estimates stems from the fact that, as it was shown in \cite{JudNem2018,juditsky2020polyhedral,PUP},
in the Gaussian case ($\xi\sim\cN(0,\sigma^2I_m)$), linear and polyhedral estimates with properly designed efficiently
computable contrast matrices are near-minimax optimal in terms of their risks over a rather general class of loss functions
and signal sets which we call {\em ellitopes}.
\footnote{Exact definitions of these sets are reproduced in the main body of the paper.
For the time being, it suffices to point out an instructive example: a bounded
intersections of finitely many sets of the form $\{x:\|Px\|_p\leq1\}$, $p\geq2$, is an ellitope.}
In this paper, our goal is to investigate numerical algorithms for design of near-optimal linear and polyhedral estimates. Specifically, we aim at developing numerical routines for efficient computation of contrast matrices $H$, the principal parameters of the estimates of both types.

As it was shown in \cite{JudNem2018,juditsky2020polyhedral,PUP}, given the problem data---matrices $A,\,B$, the signal ellitope $\cX$, and the {\em co-ellitopic norm} $\|\cdot\|$ in which estimation error is measured, computing the contrast matrices of linear and polyhedral estimates amounts to solving a well-structured convex optimization problem with linear objective and linear matrix inequality constraints. State-of-the-art optimization software, e.g., CVX \cite{cvx2014} which relies upon Interior Point Semidefinite Programming (SDP) algorithms may be used to compute high-accuracy solutions to these problems. However, the structure of the optimization problems in question (presence of ``dense'' matrix arguments) results in prohibitively long processing times by IPM algorithms already for rather moderate problem dimensions (signal dimension $n$ in the range of few dozens). In this paper we discuss an alternative approach to solving the problem of designing linear and polyhedral estimates which rely upon first order algorithm, namely, the Composite Truncated Level method (CTL) of the bundle family. In particular, we show how matrix arguments can be eliminated from the contrast optimization problem and how the problem can be cast in  the form amenable for first order algorithms.
\par
The paper is organized as follows. In Section \ref{sec:linp} we present the precise setting of the estimation problem and define optimization problems underlying the contrast design for linear and polyhedral estimates. To set up the first order optimization algorithm, we demonstrate how the problem of contrast computation for linear estimate can be reduced to that for the polyhedral one, and then explain in Section \ref{fosolve} how the latter problem can be rewritten in the form not involving matrix arguments and convenient for solving using first order algorithms. A small simulation study presented in Section \ref{sec:num} illustrates numerical performance of the proposed algorithms.
We present the details of the Composite Truncated Level algorithm in Section \ref{CTBL} of the appendix. Proofs of technical statements are put to Section \ref{sec:proofs}.

\section{Linear and polyhedral estimates}\label{sec:linp}
\subsection{Estimation problem}\label{linpolprob}
Consider the problem of recovering linear image $w=Bx$  of unknown signal $x\in\bR^n$  from noisy observation
\beq\label{obs}
\omega=Ax+\sigma\xi\\
\eeq
where $B\in \bR^{\nu\times n}$ and $A\in\bR^{m\times n}$ are given matrices and $x$ is known to belong to a given signal set $\cX$. Throughout the paper $\xi$ is $(0,I_m)$-sub-Gaussian, denoted $\xi\sim \SG(0,I_m)$, i.e.,
for all $t\in \bR^m$,
\be
\bE\left\{e^{t^T\xi}\right\}\leq \exp\left(\half {\|t\|_2^2}\right).
\ee{sg0}
Given a norm $\|\cdot\|$ on $\bR^\nu$ and reliability tolerance $\epsilon\in(0,1)$,  we quantify the performance of a candidate estimate $\widehat{w}(\cdot)$ by its {\sl $\epsilon$-risk}
\be
\Risk_{\epsilon}[\widehat{x}|\cX]&=\sup_{x\in \cX}\inf_\rho\left\{\rho: \Prob_{\xi}\left\{\|\widehat{w}(Ax+\xi)-Bx\|>\rho\right\}\leq\epsilon\right\}.
\ee{erisk}
We assume from now on that the signal set $\cX$ and the {\sl polar} $\cB_*$ of the unit ball of $\|\cdot\|$ are basic ellitopes (see, e.g., \cite{l2estimation} and \cite[Section 4.2]{PUP}).
Specifically, we set
\beq\label{cXcB}
\begin{array}{rcl}
\cX&=&\{x\in\bR^n:\exists t\in \cT: x^TT_kx\leq t_k,k\leq K\},\\
\cB_*&=&\{y\in\bR^\nu:\exists s\in\cS: y^TS_\ell y\leq s_\ell,\ell\leq L\}.\\
\end{array}
\eeq
Here $T_k\succeq0$ with $\sum_kT_k\succ0$ (respectively, $S_\ell\succeq0$ with $\sum_\ell S_\ell\succ0$,
and $\cT\subset\bR^K_+$ (respectively, $\cS\subset\bR^L_+$) is a convex compact set which is monotone (i.e., $0\leq t\leq t'\in\cT$ $\Rightarrow t'\in\cT$) and possesses a nonempty interior.
We refer to $K$ (respectively, to $L$) as {\em ellitopic dimension} of $\cX$ (respectively, of $\cB_*$).

Every basic ellitope is a convex compact set with nonempty interior which is symmetric w.r.t. the origin.
``Standard'' examples of basic ellitopes are:
\begin{itemize}
\item A bounded intersection $\cX$ of $K$ centered at the origin ellipsoids/elli\-p\-tic cylinders $\{x\in \bR^n:x^T T_kx\leq1\}$ [$T_k\succeq0$]:
\[
\cX=\{x\in \bR^n:\exists t\in\cT:=[0,1]^K: x^T T_kx\leq t_k,\,k\leq K\}
\]
In particular, the unit box $\{x\in \bR^n:\|x\|_\infty\leq1\}$ is a basic ellitope.
\item A $\|\cdot\|_p$-ball in $\bR^n$ with $p\in[2,\infty]$:
\[
\{x\in\bR^n:\|x\|_p\leq1\} =\big\{x:\exists t\in\cT=\{t\in\bR^n_+,\|t\|_{p/2}\leq 1\}:\underbrace{ x_k^2}_{x^T T_k x}\leq t_k,\,k\leq K\big\}.
\]
\end{itemize}
\subsection{The estimates}\label{sec:linpo}
The interest of ellitopes in the present context is motivated by the results of \cite{l2estimation,juditsky2020polyhedral,PUP} which state that, in the situation in question, minimax-optimal, within logarithmic in $K+L$ factor, estimates can be found among {\sl linear} and {\sl polyhedral} ones.
\subsubsection{Linear estimate}
Linear estimate is specified by an $m\times \nu$ contrast matrix $H$ according to
\[
\widehat{w}_H(\omega)=H^T\omega.
\]
Let
\begin{align}\label{rdel}
\mr_{\varkappa}(H):=\min_{\lambda,\mu,\Theta}
\bigg\{&\phi_\cT(\mu)+\phi_{\cS}(\lambda)+\sigma^2\varkappa^2 \Tr(\Theta):\,\lambda \geq 0,\mu\geq 0\\
&\left.\left[
\begin{array}{c|c|c}
\sum_\ell\lambda_\ell S_\ell&\half (B-H^TA)&\half {H^T}
\cr\hline
\half (B-H^TA)^T&\sum_k\mu_kT_k&
\cr\hline
\half {H}&&\Theta
\end{array}\right]\succeq 0\right\}\nonumber
\end{align}
where for $\cG\subset \bR^p$
\[
\phi_\cG(z)=\sup_{g\in \cG}z^Tg:\,\bR^p\to\bR\cup\{+\infty\}
\]
is the support function of $\cG$.
\begin{proposition}{(cf. \cite[Proposition 4.14]{PUP})}\label{cor:VW}
Let $\wh w^H_{\mathrm{lin}}(\omega)=H^T\omega$ with some $H\in \bR^{m\times \nu}$.
\item[(i)] Then
\[\sup_{x\in \cX}\bE\left\{\|\widehat{w}(Ax+\xi)-Bx\|\right\}\leq \mr_{1}(H)
\]
\item[(ii)] Furthermore, let
\[\varkappa=1+\sqrt{2\ln[\epsilon^{-1}]}.\]
Then
\be\Risk_{\epsilon}[\wh{w}^H_{\mathrm{lin}}|\cX]\leq \mr_{\varkappa}[H].
\ee{lrrb1}
Furthermore, function $\mr_{\varkappa}[H]$ is a convex, continuous and coercive function of the contrast matrix, and can be efficiently minimized w.r.t. $H$.
\end{proposition}
\paragraph{Remarks.}\begin{itemize}
\item As a consequence of the statement $(i)$ of the proposition, the optimal value $\mr^*$ of the (clearly solvable) convex optimization problem
\begin{align}\label{syntlin}
\mr^*:=\min_H\mr_{\varkappa}(H)
=\min_{{H,\lambda,\mu,\Theta}}
\bigg\{&\phi_\cT(\mu)+\phi_{\cS}(\lambda)+\sigma^2\Tr(\Theta):\,\lambda \geq 0,\mu\geq 0\\
&\left.\left[
\begin{array}{c|c|c}
\sum_\ell\lambda_\ell S_\ell&\half (B-H^TA)&\half {H^T}
\cr\hline
\half (B-H^TA)^T&\sum_k\mu_kT_k&
\cr\hline
\half {H}&&\Theta
\end{array}\right]\succeq 0\right\}\nonumber
\end{align}
is an upper bound on the expected risk
\[
\Risk[\widehat{w}^{*}_\tLin|\cX]=\sup_{x\in \cX}\bE\left\{\|\widehat{w}^{*}_\tLin(Ax+\xi)-Bx\|\right\}
\] of the estimate $\widehat{w}^{*}_\tLin(\omega)=H_*^T\omega$ yielded by the $H$-component of an optimal solution to the problem. {Note that we do not need to assume that the noise $\xi$ is sub-Gaussian for the above bound to hold: it would suffice to suppose that $\bE\{\xi\xi^T\}\preceq \sigma^2 I$.}
Moreover (cf. \cite[Proposition 4.16]{PUP}), the value $\mr^*$ is within moderate factor of the minimax-optimal $\|\cdot\|$-risk
\[
\RiskOpt[\cX]=\inf_{\widehat{w}(\cdot)}\Risk[\widehat{w}|\cX]
\]
(here $\inf$ is taken over all estimates, linear and nonlinear alike):
\beq
\mr^*\leq O(1)\sqrt{\ln(K+1)\ln(L+1)}\RiskOpt[\cX]
\eeq
(from now on, $O(1)$'s stands for an absolute constant).
\item Observe that the $\epsilon$-risk of the estimate $\wh w^*_\lin$ may be evaluated using the statement $(ii)$ of Proposition \ref{cor:VW}.
\begin{corollary}\label{col:dev1}
Let $\epsilon\in(0,1)$. The $\epsilon$-risk of the estimate  $\widehat{w}^{*}_\tLin(\omega)$ satisfies
\[
\Risk_{\epsilon}[\wh{w}^H_{\mathrm{lin}}|\cX]\leq \left(1+\sqrt{2\ln[\epsilon^{-1}]}\right)\mr^*.
\]
\end{corollary}
\item The bound $\ov \psi_\epsilon(\Theta)=\varkappa^2\Tr(\Theta)$ for the $1-\epsilon$ quantile of the quadratic form $\xi^T\Theta\xi$ in the expression \rf{rdel} can be replaced by the following tighter but harder to process bounds (cf, e.g., \cite[Proposition A.1]{bekri2023robust})
\begin{align}
\psi_\epsilon(\Theta)&:=\min_{\alpha}\big\{-\tfrac{\alpha}{2}\log\Det(I-2\alpha^{-1} \Theta)+\alpha\ln[\epsilon^{-1}],\;\alpha\geq 2\lambda_{\max}(\Theta)\big\}\nn
\leq \psi'_{\epsilon}(\Theta)&:=\min_{\alpha}\big\{\Tr(\Theta)+
\Tr\left(\Theta(\alpha I-2\Theta)^{-1}\Theta\right)+\alpha\ln(\epsilon^{-1}):\,\alpha\geq 2\lambda_{\max} (\Theta)\big\}
\nn\leq
\wt\psi_\epsilon(\Theta)&:=\Tr(\Theta)+2\|\Theta\|_{\mathrm{Fro}}\sqrt{\ln[\epsilon^{-1}]}+2\lambda_{\max}(\Theta)\ln[\epsilon^{-1}]\leq \ov\psi_\epsilon(\Theta).\label{barpsi}
\end{align}
 \end{itemize}
\subsubsection{Polyhedral estimate}
Polyhedral estimate  is specified by  $m\times m$ contrast matrix $H$ satisfying
\beq\label{sigmachi}
\sigma\chi\|\Col_j[H]\|_2\leq1,\;1\leq j\leq m,\eeq
with \[\chi=[\chi(\epsilon/m)=]\sqrt{2\ln[2\epsilon/m]};
\] in other words, $H$ is normalized by the requirement that
\beq
\Prob_{\xi\sim \SG(0,I)}\{\sigma \|H^T\xi\|_\infty>1\}\leq\epsilon.
\eeq
The associated with $H$ polyhedral estimate $\widehat{w}_H(\omega)$ is
\beq
\widehat{w}_H(\omega)=B\overline{x}(\omega),\,\,\overline{x}(\omega)=\argmin_x\left\{\|H^T(\omega-Ax)\|_\infty:\,x\in\cX.\right\}
\eeq
\def\tPoly{\mathrm{poly}}
It is shown in \cite[Setion 5.1.5]{PUP} that the $\epsilon$-risk of $\widehat{w}_H$ is upper-bounded by the quantity
\begin{align}
\min_{\lambda,\mu,\upsilon}\Big\{&2\Big[\phi_\cS(\lambda)+\phi_\cT(\mu)+\sum_j\upsilon_j\Big]:\,
\lambda\geq0,\mu\geq0,\upsilon\geq0\tag{$P[H]$}\\
&\left.\left[\begin{array}{c|c}\sum_\ell\lambda_\ell S_\ell&\half B\cr\hline \half B^T&A^TH\Diag\{\upsilon\}H^TA+\sum_k\mu_kT_k\\
\end{array}\right]\succeq0\right\}\nonumber
\end{align} A presumably good synthesis of the contrast $H$ is yielded by feasible solutions to the convex optimization problem
\begin{align}\label{syntpol}
\mp^*_\chi=\min_\Theta\Big\{\mp_\chi[\Theta]:=\min_{\lambda,\mu,\Theta}\Big\{&2\big[\phi_\cS(\lambda)+\phi_\cT(\mu)+\sigma^2\chi^2(\epsilon/m)\Tr(\Theta)\big]:\,
\lambda\geq0,\,\mu\geq0,\,\Theta\succeq0\\
&\qquad\left.\left[\begin{array}{c|c}\sum_\ell\lambda_\ell S_\ell&\half B\cr\hline \half B^T&A^T\Theta A+\sum_k\mu_kT_k\cr\end{array}\right]\succeq0\right\}\nonumber
\end{align}
Given a feasible solution $(\lambda,\mu,\Theta)$ to \rf{syntpol}, we set $H=[\sigma\chi(\epsilon/m)]^{-1}U$ where $\Theta=U\Diag\{\nu\}U^T$ is the eigenvalue decomposition of $\Theta$ and  $\upsilon=[\sigma\chi(\epsilon/m)]^2\nu$. Note that such $H$ satisfies (\ref{sigmachi}) and, moreover, $(\lambda,\mu,\upsilon,H)$ is a feasible solution to $(P[H])$ with values of respective objectives of both problems at these feasible solutions being equal to each other. As a result, the $\epsilon$-risk of the polyhedral estimate $\widehat{w}_\tPoly$ stemming, in the just explained fashion, from an optimal solution to the (clearly solvable) problem (\ref{syntpol}) is upper-bounded by $\mp^*_\chi$. As shown in
\cite[Proposition 5.10]{PUP}, the resulting polyhedral estimate is nearly minimax-optimal:
\[
\Risk_{\epsilon}[\widehat{w}_\tPoly|\cX]\leq \mp^*_\chi\leq O(1)\sqrt{\ln(K+1)\ln(L+1)\ln(2m/\epsilon)}\,\RiskOpt_{\epsilon}[\cX].
\]
\subsection{From  polyhedral to  linear estimate and back}
Observe that  problems (\ref{syntlin}) and (\ref{syntpol}) responsible for the design of nearly minimax-optimal under the circumstances linear and polyhedral estimates are well structured convex problems. State-of-the-art Interior Point Semidefinite Programming (SDP) algorithms may be use to compute high-accuracy solutions to these problems   in a wide range of geometries of $\cT$ and $\cS$. However, the presence of matrix variables $H$ and $\Theta$ results in large design dimensions of the SDP's to be solved and make prohibitively time consuming processing problem instances of sizes $m,n$ in the range of hundreds.
The first goal of this paper is to show that matrix variables may be eliminated from (\ref{syntlin}) and (\ref{syntpol}) allowing for processing by dedicated First Order algorithms, resulting in significant extension of the ranges of problem sizes amenable for numerical processing.

Our first observation is that problems (\ref{syntlin}) and (\ref{syntpol}) are ``nearly reducible'' to each other. Indeed, let $(\lambda,\mu,H,\Theta)$ be feasible to \rf{syntlin}. We clearly have $\Theta\succeq0$. Let
$G=\left[\begin{array}{c|c|c}I&&\cr\hline&I&A^T\end{array}\right]$. By multiplying the semidefinite constraint of \rf{syntlin} by $G$ on the left and $G^T$ on the right we conclude that $(\lambda,\mu,\Theta)$ is a feasible solution to \rf{syntpol} with the corresponding objective value
\[2[\phi_\cS(\lambda)+\phi_\cT(\mu)+\sigma^2\chi^2\Tr(\Theta)]=2\mr_{\chi}[H].
\]
The converse is also true.
\begin{lemma}\label{lem:converse}
Let $(\lambda, \mu,\Theta)$ be a feasible solution to \rf{syntpol}. Then it can be augmented to the feasible solution $(2\lambda,\mu,H,\Theta)$ of \rf{syntlin} with the corresponding objective value
\[2\phi_\cS(\lambda)+\phi_\cT(\mu)+\sigma^2\varkappa^2\Tr(\Theta)\leq \mp_\varkappa[\Theta].\]
\end{lemma}
\def\mT{{\mathfrak{T}}}
\section{Designing polyhedral estimates by a First Order method}\label{fosolve}
According to the results from the previous section, when speaking about numerical design of linear and polyhedral estimates, we can focus solely on solving problem (\ref{syntpol}).\footnote{Strictly speaking, this is so if we assume that when looking for a linear estimate, we are ready to tolerate a moderate constant factor (namely, 2) in the risk bound of the resulting estimate.}  Next, projecting, if necessary, the observation $\omega$ onto the image space of $A$, we can assume w.l.o.g. that $m\leq n$ and the image space of $A$ is the entire $\bR^m$. In fact, we make here a stronger assumption:\footnote{We briefly describe the ``conversion'' of \rf{syntpol} into the form amenable for First Order algorithms in the case of singular $A$ in Section \ref{sec:app_sing}.}
\begin{quote}
{\bf AssO}: Matrix $A\in \bR^{n\times n}$ in (\ref{obs}) is  nonsingular.
\end{quote}
Under this assumption, we can carry out partial minimization in $\Theta$ in \rf{syntpol}. Specifically, it is immediately seen that \rf{syntpol} is equivalent to the optimization problem
\begin{align}\label{syntpolnew}
\min\limits_{\lambda,\mu,\Theta}\big\{&\phi_\cS(\lambda)+\phi_\cT(\mu)+\sigma^2\chi^2\Tr(\Theta):\,\lambda>0,\,\mu\geq0,\,\Theta\succeq0\\
&\qquad\qquad\Theta\succeq\underbrace{A^{-T}\left[\four B^T\left[{\sum}_\ell\lambda_\ell S_\ell\right]^{-1}B-\sum_k\mu_kT_k\right]A^{-1}}_{\mT(\lambda,\mu)}\bigg\}\nonumber
\end{align}
In the latter problem partial minimization in $\Theta$ is as follows: given $\lambda>0$ and $\mu\geq0$ we compute the eigenvalue decomposition
$$
\mT(\lambda,\mu)=U\Diag\{\upsilon\}U^T
$$
of $\mT(\lambda,\mu)$. The best in terms of the objective of (\ref{syntpolnew}) choice of $\Theta$ given $\lambda$, $\mu$  is
$$
\Theta=U\Diag\{\upsilon^+\}U^T\eqno{[\alpha^+=\max[\alpha,0]]}.
$$
Therefore, \rf{syntpolnew} reduces to
\begin{align}\label{syntpolfin}
\min_{\lambda>0,\mu\geq0} \Big\{\Upsilon(\lambda,\mu)&:=\phi_\cS(\lambda)+\phi_\cT(\mu)+\sigma^2\chi^2\sum_i\lambda^+_i(\mT(\lambda,\mu))\Big\}\\
\mT(\lambda,\mu)&=A^{-T}\left[\four B^T\Big[\sum_\ell\lambda_\ell S_\ell\Big]^{-1}B-\sum_k\mu_kT_k\right]A^{-1}\nonumber
\end{align}
where $\lambda_1(Q)\geq\lambda_2(Q)\geq...\geq\lambda_p(Q)$ are the eigenvalues of  symmetric $p\times p$ matrix $Q$. and $\lambda_i^+(Q)=\max[\lambda_i(Q),0]$.
\subsection{Setting up Composite
Truncated Level algorithm}
We intend to solve the problem of interest \rf{syntpol} by applying to (\ref{syntpolfin}) a First Order algorithm---Composite
Truncated Level algorithm {(CTL)}. Detailed description of the method is presented in Section \ref{CTBL}. {CTL} is aimed at solving convex optimization problems of the form
\beq\label{eqprob}
\Opt=\min\limits_{x\in X}\{ \phi(x):=\psi(x)+ f(x)\}
\end{equation}
where $X\subset\bR^N$ is a nonempty bounded and closed convex set, and $\psi(\cdot)$ and $f(\cdot)$ are Lipschitz continuous convex functions on $X$ with ``simple'' $X$ and $\psi$ (for details, see Section \ref{CTBL}). Note that in order to reduce the problem of interest \rf{syntpolfin} to the form \rf{eqprob}, it suffices to set
\begin{align}\label{setup}
X=\bigg\{x&=[\lambda;\mu]\in\bR^{L+K}_+,\,\lambda_\ell\geq\delta \,\forall \ell,\,\sum_\ell\lambda_\ell+\sum_k\mu_k\leq R\bigg\},\nn
\psi([\lambda;\mu])&=\phi_\cS(\lambda)+\phi_\cT(\mu),\\
f([\lambda;\mu])&=\sigma^2\chi^2\sum_i\lambda^+_i(\mT(\lambda,\mu)).\nonumber
\end{align}
When solving (\ref{eqprob}), {CTL} ``learns'' the difficult part $f(x)$ of the objective via {\sl oracle} which, given on input a query point $\overline{x}\in X$, returns a ``simple'' Lipschitz continuous convex function (model) $f_{\overline{x}}(\cdot)$ such that
$$
f_{\overline{x}}(\overline{x})=f(\overline{x})\;\; \&\;\; f_{\overline{x}}(y) \leq f(y)\quad\forall y\in X.
$$
{\bf Oracle $\cO_\varrho$.} In the situation we are interested in with the data for (\ref{eqprob}) given by (\ref{setup}), the oracle may be built as follows:
\begin{enumerate}
\item Given query point $\ov x=[\overline{\lambda};\overline{\mu}]\in X$, we compute the matrices $\overline{\Lambda}=\sum_\ell\overline{\lambda}_\ell S_\ell$ and $\overline{\mT}=\mT(\overline{\lambda},\overline{\mu})$ along with the eigenvalue decomposition
    $
    \overline{\mT}=\overline{U}\Diag\{\overline{\upsilon}\}\overline{U}^T$,
    $\overline{\upsilon}_1\geq\overline{\upsilon}_2\geq...\geq\overline{\upsilon}_n
    $
    of $\overline{\mT}$.
\item We put
$$
T(\lambda,\mu)=A^{-T}\left[\four B^T\overline{\Lambda}^{-1}\Big[2\overline{\Lambda}-{\sum}_\ell\lambda_\ell S_\ell\Big]\overline{\Lambda}^{-1}B-{\sum}_k\mu_kT_k\right]A^{-1}
$$
This function is obtained from $\mT(\lambda,\mu)$ by linearization in $\lambda$ at $\lambda=\overline{\lambda}$ and clearly $\succeq$-underestimates $\mT(\lambda,\mu)$ on $X$, while $T(\overline{\lambda}, \mu)\equiv \mT({\ov\lambda},\mu)$. Consequently,
$$\
f([\lambda;\mu])=\sigma^2\chi^2\sum_i\lambda^+_i(\mT(\lambda,\mu))\geq f_{\overline{\lambda}}(\lambda,\mu):=\sigma^2\chi^2\sum_i\lambda_i^+(T(\lambda,\mu)),
$$
the inequality being equality when $\lambda=\overline{\lambda}$.
\par
Recall that for every symmetric convex function $g$ on $\bR^n$ and every $n\times n$ symmetric matrix $D$ one has \cite[Proposition 4.2.1]{BN2001}
$$
g([D_{11};...;D_{nn}])\leq g(\lambda(D)).
$$
When specifying $g(s)=\sum_{i=1}^ns_i^+$ and denoting by $D_i(\lambda,\mu)$ the diagonal entries in the matrix $\overline{U}^TT(\lambda,\mu)\overline{U}$, $1\leq i\leq n$, we get
$$
[f([\lambda;\mu])\geq]\,f_{\overline{\lambda}}(\lambda,\mu)\geq \sum_{i=1}^nD_i^+(\lambda,\mu) \quad\forall
(\lambda>0,\mu\geq0)
$$
with both inequalities becoming equalities for $\lambda=\overline{\lambda}$ and $\mu=\overline{\mu}$. Taking into account that functions $D_i(\lambda,\mu)$ are affine, we conclude that
{\em piecewise linear function $\sum_iD_i^+(\lambda,\mu)$ underestimates $f([\lambda;\mu])$ in the domain $\lambda>0$ and is equal to $f([\lambda;\mu])$ when $(\lambda,\mu)=(\overline{\lambda},\overline{\mu})$}.
\item The above considerations justify the oracle $\cO_\varrho$ defined as follows:
$$
f_{\overline{x}=[\overline{\lambda};\overline{\mu}]}(x)=\sum_{\iota=1}^\varrho\max[\alpha_\iota(x),0],
$$
where $\varrho$, $1\leq \varrho\leq n$, is ``complexity parameter'' of the oracle, and $\alpha_\iota(x)$ are affine functions of $x=[\lambda;\mu]$ specified according to
\begin{itemize}
\item for $\iota<\varrho$, $\alpha_\iota(x)=D_\iota(x)$;
\item $\alpha_\varrho(x)=\sum_{\iota\geq\varrho}\overline{D}_\iota(x)$, $\overline{D}_\iota(x)=\left\{\begin{array}{ll}D_\iota(x),&D_\iota(\overline{x})\geq0,\\
    0,&\hbox{otherwise}.
    \end{array}\right.$
\end{itemize}
By construction, $f_{\overline{x}}(x)$ is the sum of $\varrho$ positive parts of linear forms, underestimates $f(x)$ everywhere, and coincides with $f(x)$ when $x=\overline{x}$.
 \end{enumerate}
 \subsection{Numerical illustration}\label{sec:num}
Consider the situation in which $\|\cdot\|$ is $\|\cdot\|_2$. In this case problem (\ref{syntpol}) reads
 \beq\label{prob2ini}
 \mp_\chi^*=\min\limits_{\lambda,\mu,\Theta}\left\{2[\lambda+\phi_\cT(\mu)+
\sigma^2\chi^2\Tr(\Theta)]:
 \begin{array}{ll}
 \lambda\geq0,\mu\geq0,\Theta\succeq0\\
 \left[\begin{array}{c|c}\lambda I_\nu&\half B\cr\hline \half B^T&A^T\Theta A+\sum_k\mu_k T_k\\ \end{array}\right]\succeq0\\
 \end{array}\right\}
 \eeq
 Observe that scaling a feasible solution $(\lambda,\mu,\Theta)$ to the problem according to $(\lambda,\mu,\Theta)\mapsto(s\lambda,s^{-1}\mu,s^{-1}\Theta)$ with $s>0$ preserves feasibility; the best in terms of the objective
 scaling of $(\lambda,\mu,\Theta)$ corresponds to $s=\sqrt{[\phi_\cT(\mu)+\sigma^2\chi^2\Tr(\Theta)]/\lambda}$ and results in the value of the objective $4\sqrt{\lambda[\phi_\cT(\mu)+\sigma^2\chi^2\Tr(\Theta)]}$, As a result, we can eliminate the variable $\lambda$, thus arriving at the problem
 {\small\beq\label{prob2}
 \ov\mp_\chi^*=\min\limits_{\overline{\mu},\overline{\Theta}}\left\{F(\overline{\mu},\overline{\Theta}):=\phi_\cT(\overline{\mu})+\sigma^2\chi^2\Tr(\overline{\Theta}):
 \begin{array}{ll}
 \overline{\mu}\geq0,\overline{\Theta}\succeq0\\
 \left[\begin{array}{c|c}I_\nu&\half B\cr\hline \half B^T&A^T\overline{\Theta} A+\sum_k\overline{\mu}_k T_k\\
 \end{array}\right]\succeq0\\
 \end{array}\right\}
 \eeq}\noindent
 A feasible (an optimal) solution $\overline{\mu},\overline{\Theta}$ to (\ref{prob2}) gives rise to feasible (resp., optimal)  solution $\lambda
 =\sqrt{F(\overline{\mu},\overline{\Theta})}$, $\mu=\overline{\mu}/\lambda$, $\Theta=\overline{\Theta}/\lambda$ to (\ref{prob2ini}) with the value of the objective equal to $4\sqrt{F(\overline{\mu},\overline{\Theta})}$; in particular,
 $$
 \mp_\chi^*=4\sqrt{\ov\mp_\chi^*}.
 $$
 In our experiments, we used $B=I_n$, $A\in\bR^{n\times n}$ with i.i.d. entries drawn from $\cN(0,1)$, we put
$\sigma=0.1$, $\epsilon=0.05,$ and used  ellitopic signal set
 $$
 \cX=\{x\in\bR^n:\sum_{i\in I_k}i^\alpha x_i^2\leq 1,\,1\leq i\leq K\}
 $$
where $I_1,...,I_K$ are consecutive segments of the range $1\leq i\leq n$ of cardinality $n/K$ each; we put $\alpha=1$.
\par
Under the circumstances, problem {\rf{prob2}}  reads
\beq\label{prob2us}
 \ov\mp_\chi^*=\min\limits_{\overline{\mu},\overline{\Theta}}\left\{\sum_{k=1}^K\overline{\mu}_k+
\underbrace{\sigma^2\chi^2(\epsilon/n)}_{\gamma}
\Tr(\overline{\Theta}):
 \begin{array}{ll}
 \overline{\mu}\geq0,\overline{\Theta}\succeq0\\
 \left[\begin{array}{c|c}I_\nu&\half I_n\cr\hline \half I_n&A^T\overline{\Theta} A+D[\overline{\mu}]\\
 \end{array}\right]\succeq0\\
 \end{array}\right\}
\eeq
where $D[\overline{\mu}]$ is diagonal $n\times n$ matrix with the $i$-th diagonal entry equal to $\mu_k$ when $i\in I_k$.  \par
After eliminating $\Theta$  by partial minimization (\ref{prob2us}) becomes
\beq\label{prob2fin}
\min_{\overline{\mu}}\left\{\sum_k\overline{\mu}_k+\gamma\sum_\iota\lambda_\iota^+(A^{-T}\left[I_n-D[\overline{\mu}]\right]A^{-1}): 0\leq\overline{\mu},\sum_k\overline{\mu}_k\leq R\right\}
\eeq (we have imposed a large enough upper bound on $\sum_k\overline{\mu}_k$ to make the optimization domain bounded).
The resulting problem was processed by the {CTL} algorithm %
utilizing oracle $\cO_\varrho$. $\varrho$  was the first of the two
control parameters used the experiments; the second parameter was the maximum cardinality $\tau$ of bundle allowed for {CTL}.\footnote{For description of {CTL} and related entities, see Section \ref{CTBL}.}
We used ``$\ell_1/\ell_2$ proximal setup,'' \cite[Section 5.3.3.3.]{LMCO}, in which
$$
\|\cdot\|=\|\cdot\|_1,\;\;\omega(\overline{\mu})=\kappa_n\|\mu\|_{p_n}^2,\;\;p_n=1+1/\ln n,
$$
where $\kappa_n$
is an easy to compute constant ensuring strong convexity, modulus 1, of $\omega(\cdot)$ w.r.t. $\|\cdot\|_1$.
The {CTL} parameters $\lambda_\ell$, $\overline{\theta}$, $\underline{\theta}$ were set to $1/2$, and the auxiliary problems (steps 4.2-4.3) were
processed by Interior Point solver {\tt Mosek} invoked via  {\tt CVX}, see \cite{cvx2014}.
When solving (\ref{prob2fin}), computations were terminated when the best found so far value of the objective were within the factor 1.1 of the generated by the method lower bound on the optimal value.
\par We report on results of our experiments in Tables \ref{Tabeq1} and \ref{Tabeq2}. To put these results in proper perspective, note that solving (\ref{prob2us}) by state-of-the-art {\tt Mosek} Interior Point solver takes 35 sec when $n=64,\,K=8$ and 1785 sec when $n=128,\,K=16$; as applied to the same problem with $n=256,\,K=64$, {\tt Mosek} runs out of memory.
\begin{table}\begin{center}
\begin{tabular}{|c|cccccc|}
\hline
$n$&64&128&256&512&1024&1024\\
\hline
$K$&8&16&32&64&128&1024\\
\hline
{calls}&23&8&10&28&24&26\\
\hline
{phases}&11&4&5&14&11&11\\
\hline
\hbox{CPU, sec} &7&2&4&20&96&781\\
\hline
$(0.05$-risk&3.215&4.283&4.519&4.169&4.216&7.868\\
\hline
$\|\cdot\|_2$-risk&1.709&2.237&2.389&2.209&2.231&4.271\\
\hline
\end{tabular}
\caption{\label{Tabeq1} Solving (\ref{prob2fin}) by {CTL} with $\varrho=10,\tau=10$.}
\end{center}
\end{table}

\begin{table}\begin{center}
\begin{tabular}{|c|cc|}
\cline{2-3}
\multicolumn{1}{c|}{}&$\varrho$=1&$\varrho$=10\\
\hline
$\tau$=1&50/11/275&31/12/136\\
\hline
$\tau$=10&26/11/109&24/11/96
\\
\hline
\end{tabular}
\caption{\label{Tabeq2} Performance of {CTL}  vs. $\varrho$ and $\tau$, problem (\ref{prob2fin}) with $n=1024,K=128$.}
Data in cells: \#~of calls/\#~of phases/CPU time, sec.
\end{center}\end{table}

\paragraph{Acknowledgement.}
This work was supported by Multidisciplinary Institute in Artificial intelligence MIAI {@} Grenoble Alpes (ANR-19-P3IA-0003).

\appendix
\section{{CTL}---Composite Truncated Level algorithm}\label{CTBL}
\subsection{Situation and goal.} {CTL} is a First Order method for solving optimization problems
$$
\Opt=\min\limits_{x\in X} \phi(x):=\psi(x)+ f(x),\eqno{(\ref{eqprob})}
$$
where
\begin{itemize}
\item $X\subset\bR^n$ is nonempty, convex, closed, and bounded
\item $\psi:X\to\bR$ and $f:X\to\bR$ are Lipschitz continuous convex functions.
\end{itemize}
Our assumptions are as follows:
\par\noindent
{\bf Ass1}: $X$ is equipped with a {\sl proximal setup} composed of a norm $\|\cdot\|$ on $\bR^n$
and  a {\sl distance-generating function} $\omega:X\to\bR$ which is continuously differentiable
and strongly convex on $X$ with convexity modulus 1 w.r.t. $\|\cdot\|$:
$$
\langle \nabla\omega(x)-\nabla\omega(y),x-y\rangle \geq \|x-y\|^2\,\;\;\forall x,y\in X.
$$
Proximal setup induces {\sl Bregman distance $V_x(y)$}  on $X$ and  {\sl Bregman diameter $\Omega$} of $X$:
\bse
V_x(y)&=&\omega(y)-[\omega(x)-\langle \nabla\omega(x),y-x\rangle]\geq {1\over 2}\|x-y\|^2,\,x,y\in X\\
\Omega&=&\left[2\max\limits_{x,y\in X}V_x(y)\right]^{1/2}\geq\max\limits_{x,y\in X}\|x-y\|.
\ese
\par\noindent
{\bf Ass2}: We have at our disposal {\sl oracle $\cO$} which, given on input a query point $x\in X$, returns a {\sl piece}---a Lipschitz continuous on $X$ convex function
$$
\phi_x(\cdot)=\psi(\cdot)+f_x(\cdot):X\to \bR
$$
where $f_x(\cdot)$ belongs to some family $\cF$ of ``simple'' Lipschitz continuous convex functions on $X$. We suppose that
\[
\forall (x,y\in X):\; f_x(y)\leq f(y) \;\; \&\;\; f_x(x)=f(x)
\]
(in particular, $\phi(x)=\phi_x(x)$)
and that functions $\phi_x(\cdot):X\to\bR$ are uniformly in $x\in X$ Lipschitz continuous on $X$:
\beq\label{eqlip}
|\phi_x(u)-\phi_x(v)|\leq L_\phi \|u-v\|\,
\,\forall (u,v\in X)
\eeq
for some $L_\phi<\infty$.
\par
The simplest example of such oracle is that of family $\cF$ comprised of affine functions of $\bR^n$, and $f_x(y)=f(x)+\langle f'(x),y-x\rangle$
where $f'(x)$ is a subgradient of $f$ at $x$ (``first order oracle''). In a  less trivial example, $f(x)$ is the largest eigenvalue of a symmetric
matrix $S(x)$ which affinely depends on $x$, while $f_x(y)=\max_{i\leq k}e_i^TS(y)e_i$, where $k\leq n$ is fixed and $e_1,...,e_k$ are the $k$ leading eigenvectors of $S(x)$.
\par\noindent
{\bf Ass3}: We assume that for some positive integer $\tau$ we are able to solve efficiently problems of the form
$$
\min\limits_y\left\{\psi(y)+\max\limits_{\iota\leq \tau}f_{x_\iota}(y):y\in X,\alpha(y)\leq 0\right\}
$$
and
$$
\min\limits_y\left\{V_x(y):  y\in X,\psi(y)+\max\limits_{\iota\leq\tau}f_{x_\iota}(y)\leq\ell,\alpha(y)\leq 0\right\}
$$
where $\alpha(\cdot)$ is an affine function.
\par
We remark that the algorithm to follow is the composite version of Non-Euclidean Restricted Memory Level method \cite{NERML}
operating with $\psi\equiv0$ and the family of affine functions in the role of $\cF$; the Euclidean version of the latter algorithm  is the minimization version
of the Proximal level bundle method from
\cite{Kiwiel}.
\subsubsection{The algorithm}\label{CTBLalg}
The description of {CTL} is as follows.\par
{\bf 0.} Control parameters of the algorithms are $\lambda_\ell\in(0,1),\overline{\theta}\in(0.1),\underline{\theta}\in(0,1)$.
\par
{\bf 1.} At the beginning of an iteration of {CTL}, we have at our disposal
\\
\indent{\bf 1.1)} {\sl upper bound} $\overline{\phi}$ on $\Opt$---the best (the smallest) of the values of $\phi$ observed at the query points processed so far. These upper bounds do not increase as the iteration count grows. The query point $\widehat{y}$ with $\overline{\phi}=\phi(\widehat{y})$ is considered as the approximate solution generated so far.\\
\indent{\bf 1.2)} {\sl lower bound} $\underline{\phi}$ on $\Opt$; these lower bounds do not decrease as the iteration number grows\\
\indent{\bf 1.3)} {\sl prox-center} $\overline{x}\in X$ and {\sl level} $\ell\in(\underline{\phi},\overline{\phi})$\\
\indent{\bf 1.4)} {\sl query point} $x\in X$\\
\indent{\bf 1.5)} {\sl bundle} -- a nonempty collection $\cB$ of $\tau_\cB\leq\tau$  pieces $\phi_{\iota}(\cdot)=\psi(\cdot)+f^\iota(\cdot)\}$, $1\leq\iota\leq\tau_{\cB}$, with $f^\iota\in\cF$; positive integer {$\tau$} is a parameter of the algorithm.\\
\noindent$\bullet$ The very first iteration is preceded by {\sl initialization} where we call the oracle at a (whatever) point $x_{\hbox{\tiny ini}}\in X$ and set
$$
\overline{\phi}=\phi(x_{\hbox{\tiny ini}}),\,\underline{\phi}=\min\limits_{x\in X}[\psi(x)+f_{x_{\hbox{\tiny ini}}}(x)],\,\cB=\{\psi(x)+f_{x_{\hbox{\tiny ini}}}(x)\}
$$
Note that since the pieces reported by the oracle underestimate $\phi(\cdot)$, we do ensure $\underline{\phi}\leq \Opt$.
\par
{\bf 2.} Iterations are split into consecutive {\sl phases}, with prox-center and level common to all iterations of a phase. For a particular phase,
\\
\indent{\bf 2.1)} the prox-center $\overline{x}$ is selected at the beginning of the first iteration of the phase and can be a whatever point of $X$;\\
\indent{\bf 2.2)} the query point of the first iteration of the phase is $x=\overline{x}$,\\
\indent{\bf 2.3)} the level $\ell$ is selected at the beginning of the very first iteration of the phase as
$$
\ell=\lambda_\ell\underline{\phi}+(1-\lambda_\ell)\overline{\phi}.
$$
$\bullet$ At the beginning of the first iteration of a phase, we set
\[
\overline{\Delta}=\overline{\phi}-\ell,\;\underline{\Delta}=\ell-\underline{\phi}, \;\Delta=\overline{\Delta}+\underline{\Delta}=\overline{\phi}-\underline{\phi}.
\]
Note that {\sl the gap} $\Delta$ of a phase upper-bounds the inaccuracy in terms of the objective of the approximate solution available at the beginning of the phase.
\par
{\bf 3.} At iterations of a phase, we maintain the relation
\beq\label{eqrel}
\phi(y)\geq\ell\;\;\mbox{for all}\;\; y\in X\;\;\mbox{such that}\;\; \langle \nabla\omega(x) -\nabla\omega(\overline{x}),y-x
\rangle <0
\eeq
{
where $\overline{x},\ell$ are the prox-center and the level of the phase, and $x$ is the query point of the iteration. Note that this relation takes place at the very first iteration of a phase, since for such an iteration $x=\overline{x}$.
\par
{\bf 4.} An iteration of a phase is organized as follows:\\
\indent{\bf 4.1)} We call the oracle at the query point $x$ of the iteration, thus getting the value of the objective $\phi(x)$ and a piece $\phi_x(\cdot)$. After $\phi(x)$ is known, we\\
\begin{itemize}\itemsep=-10pt
\item[---]update the upper bound $\overline{\phi}$ and the approximate solution $\widehat{x}$:
\[
\big(\ov{\phi},\wh x\big)=\left\{\begin{array}{ll}(\phi(x),x),&\mbox{if}\;\phi(x) < \ov\phi
\\
(\ov\phi,\wh x),&\mbox{otherwise.}\end{array}\right.
\]
\item[---]update the bundle $\cB$ by adding to it the piece $\phi_x(\cdot)$ and removing, if necessary, one of ``old'' pieces to keep the number of pieces $\phi_1,...,\phi_{\tau_\cB}$ in the resulting bundle to be at most $\tau$.\\
\end{itemize}
\indent{\bf 4.2)} If $\phi(x)-\ell\leq \overline{\theta}\,\overline{\Delta}$ (``essential progress in upper bound on $\Opt$''), we terminate the phase and pass to the next one. Otherwise we solve the auxiliary problem
    \beq\label{eqaux1}
    \widetilde{\phi}=\min_y\left\{\widehat{\phi}(y):=\max_{1\leq\iota\leq\tau_\cB}\phi_\iota(y):y\in X, \langle \nabla\omega(x)-\nabla\omega(\overline{x}),y-x\rangle\geq0\right\}
    \eeq
    (as usual, $\widetilde{\phi}=+\infty$ when the right hand side problem is infeasible), and update the lower bound $\underline{\phi}$ on $\Opt$ according to
    $$
    \underline{\phi}\mapsto \max\left[\underline{\phi},\min[\widetilde{\phi},\ell]\right]
    $$
{\bf Note:}  by (\ref{eqrel}), $\phi(y)\geq \ell $ at every point $y\in X$ which is not feasible for the optimization problem in (\ref{eqaux1}). Besides this, the pieces $\phi_i(\cdot)$ in the bundle, and, consequently, the {\sl model} $\widehat{\phi}(\cdot)$, underestimate $\phi(\cdot)$ on $X$. As a result, the quantity
$\min[\wt\phi,\ell]$, and consequently the new value of $\underline{\phi}$, indeed is a lower bound on $\Opt$.\par
If updated $\underline{\phi}$ satisfies
$$
\ell-\underline{\phi}\leq \underline{\theta}\,\underline{\Delta}
$$
(``essential progress in lower bound on $\Opt$''), we terminate the phase and pass to the next one.\\
\indent{\bf 4.3)} If the iteration in question does not result in phase change, we solve the auxiliary problem
\beq\label{eqaux2}
\min\limits_y\left\{V_{\overline{x}}(y):\;y\in X,\,\widehat{\phi}(y)\leq \ell,\, \langle \nabla\omega(x)-\nabla\omega(\overline{x}),\,y-x\rangle\geq0\right\}
\eeq
take its optimal solution, $x_+$, as the new query point, and pass to the next iteration of the phase.
\par\noindent
{\bf Note:} when we need to solve (\ref{eqaux2}), we have, by construction, {$\wt\phi\leq \ell$,} so that the problem in (\ref{eqaux2}) is feasible, a feasible solution being a minimizer in \rf{eqaux1}. Thus, the new query point $x_+$ is well defined. Besides this, from the definition of $x_+$ it follows that
$$
\langle \nabla \omega(x_+)-\nabla\omega(\overline{x}), y-x_+\rangle\geq0
$$
for every feasible solution $y$ to (\ref{eqaux2}). As a result, when $y\in X$ satisfies the relation
$$
\langle \nabla \omega(x_+)-\nabla\omega(\overline{x}), y-x_+\rangle<0,
$$
$y$ is infeasible for (\ref{eqaux2}), meaning that either $\widehat{\phi}(y)>\ell$, and in such case $\phi(y)\geq \widehat{\phi}(y)\geq\ell$, or $y$ satisfies the premise in
(\ref{eqrel}), implying that $\phi(y)\geq\ell$ by (\ref{eqrel}). We see that
(\ref{eqrel}) holds true when $x$ is replaced with $x_+$, that is, (\ref{eqrel}) is maintained during the iterations.
\subsubsection{Convergence analysis}
Observe that
\begin{quote}
(!) {\sl If a phase is finite, then the gap $\Delta_+$ of the subsequent phase does not exceed a fixed fraction $\theta \Delta$ of the gap $\Delta$ of the phase in question, where}
$$
\theta=\max\left[1-\lambda_\ell\overline{\theta},\underline{\theta}+\lambda_\ell(1-\underline{\theta})\right]\in(0,1);
$$
\end{quote}
Indeed, (!) is an immediate consequence of the phase termination rules in {\bf 4.2} combined with the facts that $\underline{\phi}$ does not decrease, and $\overline{\phi}$ does not increase as the iteration count grows.\\
The following observation is crucial:
\begin{quote}
(!!) {\sl The number of iterations at a phase with gap $\Delta$ does not exceed}
\beq\label{eqarr}
\left\lceil\left({L_\phi \Omega\over \overline{\theta}\lambda_\ell\Delta}\right)^2\right\rceil
\eeq
(here $\lceil a\rceil$ stands for the upper integer part of $a$---the smallest integer greater or equal to $a$).
\end{quote}
Indeed, let $\ell$ be the level of the phase. Assume that the phase contains more that $T\geq1$ iterations, so that the upper bound $\ov\phi$, the lower bound $\underline{\phi}$ on $\Opt$, same as the model $\phi_t(\cdot)$ generated at iteration $t$ of the phase  are well defined for $t=1,...,T$, and the query points $x_t$ are well defined for $t=1,...,T+1$. By construction,
for $1\leq t\leq T$ we have
\begin{subequations}\label{atct}
\begin{align}
\phi_t(x_t)&>\ell+\overline{\theta}\,\overline{\Delta},\label{atct:a}\\
\phi_t(x_{t+1})&\leq\ell,\label{atct:b}\\
\langle \nabla V_{\overline{x}}(x_t),x_{t+1}-x_t\rangle &\geq0\label{atct:c}
\end{align}
\end{subequations}
where $\overline{x}=x_1$ is the prox-center of the phase. Indeed, when $t\leq T$,\\
--- \rf{atct:a} holds true since otherwise the phase would be terminated at its $t$-th iteration due to essential progress in upper bound on  $\Opt$, which is not the case when $t\leq T$;
\\
--- \rf{atct:b} and \rf{atct:c} hold because, by construction of $x_{t+1}$ at a non-concluding iteration $t$ of a phase, $x_{t+1}$ minimizes continuously differentiable on $X$ convex function $V_{\overline{x}}(\cdot)$ over the set
$$X\cap\{y\in X:\phi_t(y)\leq\ell\}\cap\{y\in Y:\langle\nabla V_{\overline{x}}(x_t),y-x_t\rangle\geq0\}.$$
Now note that by construction of $\phi_t(\cdot)$ and due to Assumption Ass2, $\phi_t$ is Lipschitz continuous with constant $L_\phi $ w.r.t. $\|\cdot\|$, which combines with \rf{atct:a} and \rf{atct:b} to imply that
$$\|x_t-x_{t+1}\|> L_\phi ^{-1}\overline{\theta}\cdot\overline{\Delta}=L_\phi ^{-1}\overline{\theta}\lambda_\ell\Delta.$$
The latter relation, in turn, combines with \rf{atct:c} and with inherited from $\omega(\cdot)$ strong convexity of $V_{\overline{x}}(\cdot)$  on $X$  w.r.t. $\|\cdot\|$ to imply that
$$
V_{\overline{x}}(x_{t+1})>V_{\overline{x}}(x_t)+{[\overline{\theta}\lambda_\ell]^2\Delta^2\over 2L_\phi ^2},\,1\leq t\leq T.
$$
Taking into account that $V_{\overline{x}}(x_1)=0$ and $V_{\overline{x}}(x)\leq {1
\over 2}\Omega^2$ for all $x\in X$, we arrive at (\ref{eqarr}).
\par As an immediate consequence of (!) and (!!), we get the following efficiency estimate:
\begin{proposition}\label{propeffest}
For every $\epsilon\in(0,L_\phi \Omega)$, the overall number of {CTL} iterations until an $\epsilon$-optimal, as certified by current gap,  solution to the minimization problem is built does not exceed
$$ N(\epsilon)=C(L_\phi \Omega/\epsilon)^2
$$
with $C$ depending solely on the control parameters $\lambda_\ell$, $\overline{\theta}$, and $\underline{\theta}$.
\end{proposition}
\section{Proofs for Section \ref{sec:linp}}\label{sec:proofs}
\subsection{Proof of Proposition \ref{cor:VW}}
Note that for all $x\in \cX$, the loss $\|\wh w_H-Bx\|$ of the estimate $\wh w_H$ satisfies
\begin{align*}
\|\wh w_H-Bx\|&=\|(H^TA-B)x+\sigma H^T\xi\|=\max_{u\in \cB_*}\left\{ u^T[(H^TA-B)x+\sigma H^T\xi]\right\}\\&=
\max_{u\in \cB_*}\left\{[u;x;\sigma\xi]^T\left[\begin{array}{c|c|c} &\half (B-H^TA)&\half H^T\cr\hline\half (B-H^TA)^T&&\cr\hline\half H&&\end{array}\right][u;x;\sigma\xi]\right\}\\
&\leq \max_{u\in \cB_*,x\in \cX}\left\{u^T\left[\sum_\ell\lambda_\ell S_\ell\right]u+x^T\left[\sum_k\mu_k T_k\right]x+\sigma^2\xi^T\Theta\xi\right\}
\end{align*}\
where $\mu,\lambda\geq 0$ and $\Theta\in \bS^m$ are such that
\[\left[
\begin{array}{c|c|c}
\sum_\ell\mu_\ell S_\ell&\half (B-H^TA)&\half H
\cr\hline
\half (B-H^TA)^T&\sum_k\lambda_kT_k&
\cr\hline
\half H^T&&\Theta
\end{array}\right]\succeq 0.
\]
We conclude that for all $x\in \cX$,
\[
\|\wh w_H-Bx\|\leq \max_{s\in \cS,t\in \cT} \left\{\sum_\ell \lambda_\ell s_\ell+ \sum_k \mu_k t_k+\sigma^2\xi^T\Theta\xi\right\}=\phi_\cS(\lambda)+\phi_\cT(\mu)+\sigma^2\xi^T\Theta\xi.
\]
Due to \rf{sg0} we have $\bE\{\xi\xi^T\}\preceq\sigma^2 I$. We conclude that $\bE\{\xi^T\Theta\xi\}\leq \Tr(\Theta)$ what implies the first claim of the proposition.
To complete the proof it remains to recall the bound
\[
\Prob_\xi\left\{\xi^T\Theta\xi\geq \varkappa^2\Tr(\Theta))\right\}\leq \epsilon
\]
for deviations of the quadratic form of sub-Gaussian random vectors (cf., e.g., \cite{hsu2012tail,rudelson2013hanson,spokoiny2013sharp}).\qed
\subsection{Proof of Corollary \ref{col:dev1}}
Indeed, let $\varkappa=1+\sqrt{2\ln[\epsilon^{-1}]}$, and let $\lambda^*,\mu^*$ and $\Theta^*$ be components of an optimal solution to \rf{syntlin}.
Notice that $\lambda^*,\mu^*$ and $\Theta^*$ are feasible for \rf{rdel}. Moreover, $\ov\lambda,\ov\mu$ and $\ov\Theta$ where
\[\ov\lambda=\varkappa\lambda^*,\;\;\ov\mu=\varkappa^{-1}\mu^*,\;\; \ov\Theta=\varkappa^{-1} \Theta^*
\]
are also feasible, so, by item $(ii)$ of Proposition \ref{cor:VW}, the value
\[
\phi_\cT(\ov\mu)+\phi_{\cS}(\ov\lambda)+\sigma^2\varkappa^2 \Tr(\ov\Theta)\leq \varkappa\left(\phi_\cT(\mu^*)+\phi_{\cS}(\lambda^*)+\sigma^2 \Tr(\ov\Theta)\right)\leq \varkappa r^*
\]
upper-bounds the $\epsilon$-risk of $\wh w^*_\lin$.\qed

\subsection{Proof of Lemma \ref{lem:converse}}
Let $(\lambda,\mu,\Theta)$ be a feasible solution to \rf{syntpol} such that
\beq\label{strict}
\Lambda:=\sum_\ell\lambda_\ell S_\ell\succ0, \,\Xi :=\sum_k\mu_kT_k\succ0.
\eeq
Note that every feasible solution to \rf{syntpol} remains feasible and satisfies \rf{strict} after replacing zero entries of $\lambda_\ell$ and $\mu_k$, if any, with arbitrarily small positive entries.\par
Now, due to
$\left[\begin{array}{c|c}\Lambda&\half B\cr\hline
\half B^T&A^T\Theta A+\Xi \end{array}\right]\succeq0
$
we have
\[\left[\begin{array}{c|c}I_\nu&\half \Lambda^{-1/2}B\Xi ^{-1/2}\cr\hline
\half \Xi ^{-1/2}B^T\Lambda^{-1/2}&\Xi ^{-1/2}A^T\Theta A\Xi ^{-1/2}+I_n\end{array}\right]\succeq0.
 \]
When setting $\Theta^{1/2}A\Xi ^{-1/2}=US$ with $UU^T=I_n$ and $S\succeq0$, we get
\[\four [\Lambda^{-1/2}B\Xi ^{-1/2}]^T[\Lambda^{-1/2}B\Xi ^{-1/2}]\preceq S^2+I_n\preceq(S+I_n)^2
\]
(recall that $S\succeq 0$). We conclude that there is $Q\in \bR^{\nu\times n}$ of spectral norm $\|Q\|_{2,2}\leq 1$ such that
\[\half [\Lambda^{-1/2}B\Xi ^{-1/2}]=Q(S+I_n)
\]
and
\[
B-2\Lambda^{1/2}QS\Xi^{1/2}=2\Lambda^{1/2}Q\Xi^{1/2}.
\]
When recalling what $S$ is and setting $H=2\Theta^{1/2}UQ^T\Lambda^{1/2}$ we have
\[B-H^TA=2\Lambda^{1/2}Q\Xi^{1/2}.
\]
Due to $\|Q\|_{2,2}\leq 1$, by Schur complement lemma, now it follows that
\be
\left[\begin{array}{c|c}\Lambda&\half (B-H^TA)\cr\hline \half (B-H^TA)^T&\Xi \cr\end{array}\right]\succeq0.
\ee{p2l1}
Besides this, by construction,
\[\four H\Lambda^{-1}H^T=\Theta^{1/2}\underbrace{UQ^TQU^T}_{\preceq I_n}\Theta^{1/2}\preceq \Theta,
\]that is,
\beq\label{eq183}
\left[\begin{array}{c|c}\Lambda&\half H^T\cr\hline\half H&\Theta\cr\end{array}\right]\succeq0.
\eeq
Finally, \rf{p2l1} together with \rf{eq183} imply that matrix
\[
\left[\begin{array}{c|c|c}2\Lambda&\half (B-H^TA)&\half H^T\cr\hline \half (B-H^TA)^T&\Xi& \cr\hline
\half H&&\Theta\end{array}\right]
\]
is positive semidefinite, meaning that $2\lambda,\mu,H$ and $\Theta$ form a feasible solution to \rf{syntlin}. The corresponding objective value is
\[\phi_\cS(2\lambda)+\phi_\cT(\mu)+\sigma^2\varkappa^2\Tr(\Theta)=2\phi_\cS(\lambda)+\phi_\cT(\mu)+\sigma^2\varkappa^2\Tr(\Theta)\leq \mp_\varkappa[\Theta].\eqno{\mbox{\qed}}
\]
\section{Contrast synthesis for the polyhedral estimate }\label{sec:app_sing}
Projecting, if necessary, the observations onto the image
space of $A$, we reduce the situation to that in which this image space is the entire
$\bR^m$; that is, $m \leq n$ with positive singular values $\sigma_1,...,\sigma_m$ of A. Let us assume that $n=m+d$ for some $d\geq 0$, and let
\[A=UDV^T,\quad D=
\big[\underbrace{\Diag(\sigma_1,...,\sigma_m)}_{=:D_m},0_{m\times d}\big],\;U\in \bR^{m\times m},\;V\in\bR^{m\times n},
\]
be the (full) singular-value decomposition of $A$. The starting point of the following computation is the bound \rf{syntpol} for the $\epsilon$-risk of the polyhedral estimate: we have
\begin{align}
\mp^*_\chi&=
2\min_{\lambda,\mu, \Theta}\left\{\phi_\cS(\lambda)+\phi_\cT(\mu)
+\sigma^2\chi^2\Tr(\Theta):\,\lambda\geq 0,\mu\geq 0,\Theta\succeq 0,\right.\nn
&\qquad\qquad\qquad\qquad\left.\left[\begin{array}{c|c}\sum_\ell \lambda_\ell S_\ell
&\half B\cr\hline
\half B^T&\sum_k \mu_kT_k+A^T\Theta A\end{array}\right]\succeq 0
\right\}\nn
&=
2\min_{\lambda,\kappa, \Theta}\left\{2\phi_\cS(\lambda)+\phi_\cT(\mu)
+\sigma^2\chi^2\Tr(\Theta):\,\lambda> 0,\mu\geq 0,\Theta\succeq 0,\right.\nn
&\left.\qquad\qquad A^T\Theta A\succeq \tfrac{1}{4}B^T\Big[\sum_\ell \lambda_\ell S_\ell\Big]^{-1}B -\sum_k \mu_kT_k
\right\}\nn&=
2\min_{\lambda,\mu, \ov\Theta}\Big\{F(\lambda,\mu,\ov\Theta):=\phi_\cS(\lambda)+\phi_\cT(\mu)
+\sigma^2\chi^2\Tr(\ov\Theta):\,\lambda> 0,\mu\geq 0,\,\ov\Theta=U^T\Theta U\succeq 0,\nn
&\qquad\qquad \left[\begin{array}{c|c}D_m\ov \Theta D_m&\cr\hline &\end{array}\right]=D^T\ov\Theta D\succeq V^T\left[\tfrac{1}{4}B^T\big[\sum_\ell \lambda_\ell S_\ell\big]^{-1}B
-\sum_k \mu_kT_k\right]V
\Bigg\}\nn
&=
2\min_{\lambda,\mu, \ov\Theta}\Big\{F(\lambda,\mu,\ov\Theta):\,\lambda> 0,\mu\geq 0,\,\ov\Theta\succeq 0,\label{barthe}\\
& \left[\begin{array}{c|c}\ov \Theta &\cr\hline &\end{array}\right]\succeq \underbrace{\left[\begin{array}{c|c}D_m &\cr\hline &I\end{array}\right]^{-1}V^T\left[\tfrac{1}{4}B^T
\big[\sum_\ell \lambda_\ell S_\ell\big]^{-1}B -\sum_k \mu_kT_k\right]V\left[\begin{array}{c|c}D_m &\cr\hline &I\end{array}\right]^{-1}}_{=:C(\lambda,\mu)}.\nonumber
\Bigg\}
\end{align}
Observe that the matrix-valued function $C(\lambda,\mu)$ is $\succeq$-convex for $\lambda>0$, and is negative definite for every fixed $\lambda$  for all $\mu$ such that $\min_i\mu_i\geq \ul\mu$ large enough. On the other hand, when $Z(\lambda,\mu)\prec 0$ in the representation
\[
C(\lambda,\mu)=\left[\begin{array}{c|c}X(\lambda,\mu) &Y(\lambda,\mu)\cr\hline Y^T(\lambda,\mu) &Z(\lambda,\mu)\end{array}\right]
\]
the semidefinite constraint of \rf{barthe} is satisfied if and only if
\[\ov \Theta\succeq W(\lambda,\mu):=X(\lambda,\mu)+Y^T(\lambda,\mu)Z(\lambda,\mu)^{-1}Y(\lambda,\mu).
\]
As a result, when denoting $[M]_+$ the matrix  obtained from a symmetric matrix $M$ by replacing its eigenvalues with their positive parts in the eigenvalue decomposition of $M$, we conclude that
\begin{align*}
\mp^*_\chi&=2\min_{\lambda,\mu}\big\{\phi_\cS(\lambda)+\phi_\cT(\mu)
+\sigma^2\chi^2\Tr[W(\lambda,\mu)]_+:\,\lambda> 0,\mu\geq 0,\,Z(\lambda,\mu)\prec 0\Big\}.
\end{align*}
The bottom line is that in the situation of this section, building $\Theta$ (and thus---the near-optimal polyhedral estimate) reduces to solving convex problem of design dimension $L+M$  with (relatively) easy-to-compute objective and constraints.

\end{document}